\documentclass[12pt]{article}
\usepackage{mathrsfs}
\usepackage{stmaryrd}
\usepackage{amsfonts,amsmath,amssymb,amscd}
\usepackage{shadow}
\usepackage{graphicx,floatrow}
\usepackage{color}
\usepackage{longtable}
\usepackage{pstricks,multido}
\usepackage{hyperref}
\usepackage{qtree}
\usepackage{tabularx}
\usepackage{array}
\allowdisplaybreaks

\newtheorem{thm}{Theorem}[section]
\newtheorem{lem}[thm]{Lemma}

\newtheorem{cor}[thm]{Corollary}
\newtheorem{conj}[thm]{Conjecture}

\def\qed{\hfill \rule{4pt}{7pt}}

\def\pf{\noindent {\it{Proof.} \hskip 2pt}}

\numberwithin{equation}{section}

\begin{document}
\begin{center}
{\large\bf  On Permutations with Bounded Drop Size}
\end{center}

\begin{center}
Joanna N. Chen$^1$, William Y.C. Chen$^2$

$^{1}$Center for Combinatorics, LPMC-TJKLC\\
Nankai University\\
Tianjin 300071, P.R. China

$^2$Center for Applied Mathematics\\
Tianjin University\\
Tianjin 300072, P.R. China

$^1$joanna@cfc.nankai.edu.cn,
$^2$chenyc@tju.edu.cn

\end{center}

\begin{abstract}
The maximum drop size of a permutation $\pi$ of $[n]=\{1,2,\ldots, n\}$
is defined to be the maximum value of $i-\pi(i)$.
Chung, Claesson, Dukes and Graham
obtained
polynomials $P_k(x)$ that can be
used to determine  the number of
permutations of $[n]$ with $d$
descents and maximum drop size not larger than $k$.
Furthermore, Chung and Graham
gave combinatorial interpretations of the coefficients of $Q_k(x)=x^k P_k(x)$ and $R_{n,k}(x)=Q_k(x)(1+x+\cdots+x^k)^{n-k}$,
and raised the question of finding a bijective proof of the symmetry property of $R_{n,k}(x)$.
In this paper, we establish a bijection $\varphi$ on $A_{n,k}$, where $A_{n,k}$ is the set of permutations of  $[n]$ and maximum drop size not larger than $k$.
The map $\varphi$ remains to be a bijection between certain subsets of $A_{n,k}$. 
This provides an answer to the question of
Chung and Graham. The second result of this paper is a proof of
a conjecture of  Hyatt concerning the unimodality of polynomials
in connection with the number of signed permutations
 of $[n]$ with  $d$ type $B$ descents and the type $B$ maximum drop size not greater than $k$.
 \end{abstract}



\section{Introduction}

This paper is concerned with the study of permutations of $[n]=\{1,2,\ldots, n\}$ with $d$ descents and maximum drop size not greater than $k$. Let this number be denoted by $E^k(n,d)$.  Chung, Claesson, Dukes and Graham \cite{ChungClaesson}
found  polynomials $P_k(x)$ that can be used to determine the number
$E^k(n,d)$. They proved that the polynomials $P_k(x)$ are
 unimodal. Furthermore, Chung and Graham found combinatorial
 interpretations of the polynomials
 $Q_k(x)=x^k P_k(x)$ and $R_{n,k}(x)=Q_k(x)(1+x+\cdots+x^k)^{n-k}$, and asked for a combinatorial interpretation
 of the symmetry property of $R_{n,k}(x)$.
 The first result of this paper is to present a bijection in answer to the question of Chung and Graham.
  The second result of this paper is a proof of a conjecture
 of Hyatt \cite{Hyatt} on the unimodality of the type $B$ analogue of
 the polynomials $P_k(x)$.

Let us give an overview of  notation and terminology.
Let $S_n$ denote the set of   permutations of $[n]=\{1,2,\ldots,n\}$.
For a permutation $\pi=\pi_1 \pi_2 \cdots \pi_n$ in $S_n$, we say that a number $1\leq i \leq n-1$ is a descent of $\pi$
 if $\pi_i > \pi_{i+1}$. The descent set of $\pi \in S_n$, denoted by
 $\mathrm{Des(\pi)}$, is defined by
\begin{equation*}
\mathrm{Des(\pi)}=\{i \in [n-1]: \pi_i>\pi_{i+1}\}.
\end{equation*}
Let $\mathrm{des(\pi)}$ denote the number of descents of  $\pi \in S_n$.
We say that $\pi \in S_n$ has a drop at $i$ if $\pi_i <i$ and the drop size is meant to be $i-\pi_i$.
Define the  maximum drop size of $\pi$ by
\begin{equation*}
\mathrm{maxdrop(\pi)}=\max\{i-\pi_i\colon 1 \leq i \leq n\}.
\end{equation*}

 Chung, Claesson, Dukes and Graham \cite{ChungClaesson} obtained
 a polynomial $P_k(x)$ that can be used to
 determine the number $E^k(n,d)$ of permutations of $[n]$ with $d$ descents
and maximum drop size not larger than $k$. Let $A_{n,k}$
denote the set of permutations of $[n]$ with maximum
drop size not larger than $k$.
The $k$-maxdrop-restricted descent polynomial is defined by
\begin{equation*}
A_{n,k}(y)= \sum_{\pi \in A_{n,k}} y^{\mathrm{des(\pi)}}=\sum_{d\geq0 } E^k(n,d) y^d.
\end{equation*}
Clearly, for $k \geq n-1$,  we have $A_{n,k}=S_n$ and $A_{n,k}(y)$ becomes the Eulerian polynomial
\begin{equation*}
    A_n(y)=\sum_{\pi \in S_n} y^{\mathrm{des(\pi)}}.
\end{equation*}

Notice that the above definition of the
Eulerian polynomial differs from the definition
 given in  Stanley \cite{rp2001} by a factor of $y$.
Chung, Claesson, Dukes and Graham \cite{ChungClaesson}
obtained the following
recurrence relation of  $A_{n,k}(y)$.

\begin{thm}\label{Arecur}
For $n \geq 0$, we have
\[A_{n+k+1,k}(y)=\sum_{i=1}^{k+1} {{k+1} \choose i}(y-1)^{i-1}A_{n+k+1-i,k}(y),
\]
where  $A_{i,k}(y)=A_i(y)$ for $0 \leq i \leq k$.
\end{thm}

Using the recurrence of $A_{n,k}(y)$  in Theorem \ref{Arecur}, they
introduced the polynomials
\begin{equation}
P_k(x)=\sum_{l=0}^{k}A_{k-l}(x^{k+1})(x^{k+1}-1)^l\sum_{i=l}^{k}{i\choose l}x^{-i}
\end{equation}
to determine the numbers
$E^k(n,d)$.

\begin{thm} For $n, k\geq 0$, we have
\begin{equation}
 A_{n,k}(y)=\sum_d\beta_k((k+1)d)y^{d},\end{equation}
 where
\begin{equation}
\sum_j\beta_k(j)x^j=P_k(x)\left(\frac{1-x^{k+1}}{1-x}\right)^{n-k}.
\end{equation}
\end{thm}

In other words, $E^k(n,d)$ equals the coefficient
of $x^{(k+1)d}$ in
\begin{equation*}
P_k(x)(1+x+x^2+\dots+x^k)^{n-k}.
\end{equation*}

 Chung, Claesson, Dukes and Graham \cite{ChungClaesson} also proved that $P_k(x)$ is a unimodal polynomial. We say a sequence $s_1 s_2 \cdots s_n $ is unimodal if there exists a $1 \leq t \leq n$ such that $s_1 \leq s_2 \leq \cdots \leq s_t$ and $s_t \geq s_{t+1} \geq \cdots  \geq s_n$. A polynomial is said to be  unimodal if the sequence of its coefficients is unimodal.

Furthermore, Chung and Graham \cite{ChungGraham} found combinatorial interpretations
of  the coefficients of $Q_k(x)=x^k P_k(x)$ and $R_{n,k}(x)=Q_k(x)(1+x+\cdots+x^k)^{n-k}$.

\begin{thm}
For $n \geq 0$,
\begin{equation*}
    Q_n(x)=\sum_{0 \leq i,j \leq n} E(n+1,i ;j+1) x^{(n+1)i+j},
\end{equation*}
where $E(n,i;j)$
is the number of $\pi=\pi_1 \pi_2 \cdots \pi_n \in S_n$ such that $\mathrm{des(\pi)}=i$ and $\pi_n=j$.
\end{thm}

\begin{thm}
For $k \geq 0$,
\begin{equation*}
    R_{n,k}(x)=\sum_{0 \leq i \leq n}{\sum_{0 \leq j \leq k}} E^k(n+1,i;n+1-k+j)x^{(k+1)i+j}.
\end{equation*}
where $E^k(n,i;j)$
is the number of $\pi=\pi_1 \pi_2 \cdots \pi_n\in A_{n,k}$ such that $\mathrm{des(\pi)}=i$ and $\pi_n=j$.
\end{thm}

Notice that Chung and Graham \cite{ChungGraham} used the notation ${\left\langle\begin{matrix}n\\i \end{matrix}\right\rangle}^{j}$  for the number $E(n,i;j)$ and the notation ${\left\langle\begin{matrix}n\\i \end{matrix}\right\rangle}^{j}_{[k]}$ for the number $E^k(n,i;j)$.
 They raised the question of finding  a bijective proof of  the following symmetry property of $R_{n,k}(x)$.

\begin{thm}\label{thm-bijection}
For $n, k\geq 0$, the polynomials $R_{n,k}(x)$ are symmetric.
In other words, for  $r=(k+1)i+j$ and $r'=(k+1)i'+j'$, where $ 0 \leq i,i' \leq n-1, 0 \leq j,j' \leq k$ and
$r+r'=(n+1)k$, we have
\[
E^k(n,i;n-k+j)=E^k(n,i';n-k+j').
\]
\end{thm}

In Section \ref{sec:2}, we   construct a bijection $\varphi$ on $A_{n,k}$ by a recursive procedure. Then we prove that  $\varphi$ remains to be a bijection
between certain subsets of $A_{n,k}$. This leads to a bijective proof of Theorem \ref{thm-bijection}.

 Hyatt \cite{Hyatt} extended the notion of  the maximum drop  to the
 type $B$ case or signed permutations.
Recall that a signed permutation $\pi=\pi_1\pi_2\cdots\pi_n$
can be viewed as a permutation of $[n]$ for which
each element may be associated with a minus sign.
The descent set of a signed permutation $\pi$ is defined by
\begin{equation*}
    \mathrm{Des_{B}(\pi)}=\{i\in[0,n-1]:\pi_i>\pi_{i+1}\},
\end{equation*}
where we assume that $\pi_0=0$, see Brenti \cite{brenti}.
We denote by $\mathrm{des_B(\pi)}$ the number of type $B$ descents of $\pi \in B_n$.

Let $\mathrm{maxdrop_B(\pi)}$ denote the maximum drop size of $\pi \in B_n$ as defined by
\begin{equation*}
    \mathrm{maxdrop_B(\pi)}=\max\Big\{\max\{i-\pi_i:\pi_i>0\},\max\{i:\pi_i<0\}\Big\}.
\end{equation*}
Let $B_{n,k}$
denote the set of signed permutations of $[n]$ with maximum
drop size not larger than $k$,  and  let $E_B^k(n,d)$ denote the number of signed permutations in $B_{n,k}$ with $d$  descents.

The type $B$ $k$-maxdrop-restricted descent polynomial is defined by
\begin{equation*}
    B_{n,k}(y)=\sum_{\pi\in {B}_{n,k}}y^{\mathrm{des_B(\pi)}}=\sum_{d\geq 0} E^k_B(n,d)y^d.
\end{equation*}
When $k\geq n$, ${B}_{n,k}=B_n$ and $B_{n,k}(y)$ becomes the type $B$ Eulerian polynomial $B_n(y)$, which is defined by
\begin{equation*}
    B_n(y)=\sum_{\pi\in B_n}y^{\mathrm{des_B(\pi)}}.
\end{equation*}
 Hyatt \cite{Hyatt} showed that $B_{n,k}(y)$ satisfied the following recurrence relation.

\begin{thm}\label{Brecur}
For $n \geq 0$, we have
\[B_{n+k+1,k}(y)=\sum_{i=1}^{k+1} {{k+1} \choose i}(y-1)^{i-1}B_{n+k+1-i,k}(y),
\]
where  $B_{i,k}(y)=B_i(y)$ for $0 \leq i \leq k$.
\end{thm}

Similarly, using the above recurrence relation of $B_{n,k}(y)$, Hyatt
gave the following type $B$ analogue of the polynomials $P_k(x)$
\begin{equation}
T_k(x)=\sum_{l=0}^{k}B_{k-l}(x^{k+1})(x^{k+1}-1)^l\sum_{i=l}^{k}{i\choose l}x^{-i},
\end{equation}
which can be used to compute the numbers
$E_B^k(n,d)$.
\begin{thm}

We have
\begin{equation}
B_{n,k}(y)=\sum_d\gamma_k((k+1)d)y^{d},
\end{equation}
where
\begin{equation}
    \sum_j\gamma_k(j)x^j=T_k(x)\left(\frac{1-x^{k+1}}{1-x}\right)^{n-k}.
\end{equation}
\end{thm}

The following conjecture was posed by Hyatt \cite{Hyatt}.

\begin{conj} \label{conj-unimodal}
The polynomial $T_k(x)$ is  unimodal for any $k\geq 0$.
\end{conj}

The second result of this paper is a proof of the above conjecture,
which will be given in Section 3.

%

\section{A bijective proof of Theorem \ref{thm-bijection}}\label{sec:2}

In this section, we give a bijection $\varphi$ on $A_{n,k}$, and we prove that the map $\varphi$ remains to be a bijection
between certain subsets of $A_{n,k}$. This yields a bijective proof of Theorem \ref{thm-bijection}.

We begin with some notation.
Given $\pi \in S_n$ and $1 \leq i \leq n+1$, let $\pi \leftarrow i$ denote the permutation $\mu=\mu_1 \mu_2 \cdots \mu_{n+1}$ in
$S_{n+1}$ that is obtained from $\pi$ by adding $i$ at the end of $\pi$ and increasing the elements $i, i+1,  \ldots, n$  by $1$. For example, $3421 \leftarrow 3= 45213$.

For $n,k \geq 0$,  $0 \leq i \leq n-1$ and $0 \leq j \leq k$, let $\Gamma^k(n,i;j)$ denote the set of permutations $\pi=\pi_1\pi_2\cdots \pi_n$ in $A_{n,k}$ with $\mathrm{des(\pi)}=i$ and $\pi_n=n -k+j$, which is counted by $E^k(n,i;n-k+j)$.

Given $n$ and $k$,
we proceed to construct a map $\varphi$ from $A_{n,k}$ to
$A_{n,k}$, which can be described as  a recursive procedure.
Let  $\pi=\pi_1\pi_2\cdots \pi_n \in A_{n,k}$.
For $n=1$, we define $\varphi(1)=1$.
For $n \geq 2$,
let $i=\mathrm{des(\pi)}$ and $j=\pi_n-n+k$. Let
\begin{eqnarray}
i' &=& \Big{\lfloor} \frac{(n+1)k-(k+1)i-j}{k+1}\Big{\rfloor} \label{i'},\\[3pt]
j'&=& (n+1)k-(k+1)i-j-(k+1)i'. \label{j'}
\end{eqnarray}
Assume that $\pi'$ is the permutation of $[n-1]$
that is order isomorphic to $\pi_1\pi_2\cdots \pi_{n-1}$.
In other words, $\pi=\pi' \leftarrow \pi_n$.
Now we define $\varphi(\pi)=\varphi(\pi') \leftarrow (n-k+j')$.

For example, let $\pi=12354$, which belongs to
 $A_{5,2}$, or more precisely,
$\Gamma^2(5,1;1)$. It is easy to check that $i'=2$ and $j'=2$.
Then we have $\pi'=1234$. Iterating the above procedure, we get $\pi''=123$, $\pi'''=12$ and $\pi''''=1$.
Then we have $\varphi(\pi'''')=1$, $\varphi(\pi''')=21$, $\varphi(\pi'')=321$,
 $\varphi(\pi')=3214$ and $\varphi(\pi)=32145$.

The following theorem shows that the map $\varphi$ becomes
a bijection between certain subsets of $A_{n,k}$.

\begin{thm}\label{thmbij}
For $n,k \geq 0$,  $0 \leq i \leq n-1$ and $0 \leq j \leq k$,  the map
$\varphi$ gives
 a bijection from $\Gamma^k(n,i;j)$ to $\Gamma^k(n,i';j')$, where
$i'$ and $j'$ are given by  (\ref{i'}) and (\ref{j'}).
\end{thm}

The following lemmas are needed to prove
 Theorem \ref{thmbij}.

\begin{lem} \label{welldefined}
For $n,k \geq 0$, $0 \leq i \leq n-1$ and $0 \leq j \leq k$,
let $\pi=\pi_1 \pi_2 \cdots \pi_n$ be a permutation
in  $\Gamma^k(n,i;j)$. Then we have $\varphi (\pi)$ is a permutation
 in $\Gamma^k(n,i';j')$, where $i'$ and $j'$ are given by (\ref{i'}) and (\ref{j'}).
\end{lem}

\pf
We proceed by induction on $n$. Clearly, when
$n=1$, we have $1 \in \Gamma^k(1,0;k)$. By (\ref{i'}) and (\ref{j'}), we
have $i'=0$ and $j'=k$.  Hence $\varphi(1)=1 \in \Gamma^k(1,0;k)$.  Assume that the lemma  holds for $n-1$, where $n \geq 2$. We aim to show that it holds for $n$.

Assume that $\pi=\pi_1 \pi_2 \cdots \pi_n $ is a permutation in $ \Gamma^k(n,i;j)$, that is, $i=\mathrm{des(\pi)}$ and $j=\pi_n-n+k$.
Let $\sigma=\sigma_1 \sigma_2 \cdots \sigma_{n-1}$ be the permutation of $[n-1]$ which is order isomorphic to $\pi_1\pi_2 \cdots\pi_{n-1}$, that is, $\pi=\sigma \leftarrow \pi_n$. Then
 we have $p=p_1 p_2 \cdots p_n=\varphi(\pi)=\varphi(\sigma)\leftarrow (n-k+j')$, where $j'$ is defined by (\ref{i'}) and (\ref{j'}).

By the definition of $\varphi$, we have $p_n=n-k+j'$.
Let $\alpha=\alpha_1 \alpha_2 \cdots \alpha_{n-1}= \varphi(\sigma)$.
By the induction hypothesis, we find $\mathrm{maxdrop(\alpha)} \leq k$.
It can be seen that $\mathrm{maxdrop}(p)\leq k$ since
$p_i \geq \alpha_i$ for $1 \leq i \leq n-1$ and $n-p_n=k-j'$.

It remains to show that $\mathrm{des}(p)=i'$.
Let
\begin{eqnarray}
  a &=&\mathrm{des(\sigma)} \\[3pt]
  b &=& \sigma_{n-1}-n+1+k \\[3pt]
  a' &=& \Big{\lfloor} \frac{nk-a(k+1)-b}{k+1} \Big{\rfloor} \\[4pt]
  b' &=& nk-a(k+1)-b-a'(k+1) \label{b'}
\end{eqnarray}
Again, by the induction hypothesis, we have $a'=\mathrm{des(\alpha)}$ and $b'=\alpha_{n-1}-n+1+k$. Recall that $p=\alpha \leftarrow (n-k+j')$.
 It suffices  to prove that $i'=a'+1$ when $\alpha_{n-1} \geq p_n$ and $i'=a'$ when $\alpha_{n-1} < p_n$.
  Since $p_n=n-k+j'$ and $\alpha_{n-1}=n-1-k+b'$, we need to show that $i'=a'+1$ when $j'-b' \leq -1$ and $i'=a'$ when $j'-b' > -1$.

By the definition of $j$ and $b$, we have $0 \leq j \leq k$ and $0 \leq b \leq k$. It follows  that
\begin{equation} \label{jb}
 -k \leq j-b \leq k.
\end{equation}
Similarly, we have
\begin{equation} \label{j'b'}
 -k \leq j'-b' \leq k.
\end{equation}
By (\ref{j'}) and (\ref{b'}), we get
\begin{eqnarray}
i(k+1)+j + i'(k+1)+j' &=&(n+1)k, \label{ij} \\[3pt]
a(k+1)+b + a'(k+1)+b' &=& nk.\label{ab}
\end{eqnarray}
Here are two cases.

\noindent Case 1: $i=a$, namely, $j-b> -1$.

By (\ref{ij}) and (\ref{ab}), we find that
  \begin{equation*}
    (i'-a')(k+1)=k-(j-b)-(j'-b').
  \end{equation*}

If $j'-b' \leq -1$, by (\ref{j'b'}), we see that  $k \geq 1$.
Moreover, by (\ref{jb}) and (\ref{j'b'}) and the assumptions
$j-b>-1$ and $j'-b'\leq -1$,  we deduce that $-1< j-b \leq k$ and $-k \leq j'-b' \leq -1$. It follows that  $(i'-a')(k+1) \in [1,2k]$, where $k \geq 1$. Hence we  arrive at the assertion that $i'=a'+1$.

If $j'-b' > -1$,  by (\ref{jb}) and (\ref{j'b'}) and the assumptions $j-b>-1$, we find that $-1< j'-b' \leq k$ and $-1< j-b \leq k$. It follows that $(i'-a')(k+1) \in [-k,k]$. Hence we deduce that $i'=a'$.

\noindent Case 2: $i=a+1$, namely, $j-b \leq -1$

By (\ref{jb}) and the assumption $j-b \leq -1$, we deduce that $k \geq 1$.
It follows from (\ref{ij}) and (\ref{ab}) that
  \begin{equation*}
    (i'-a')(k+1)=-1-(j-b)-(j'-b').
  \end{equation*}

 If $j'-b' \leq -1$, we claim that $k \geq 2$.
 Assume to the contrary that $k=1$. By (\ref{jb}) and (\ref{j'b'}),  we have $j'-b'=-1$ and $j-b=-1$. Hence we get $2(i'-a')=1$, a contradiction
 to the fact that $i'-a'$ is an integer. This proves that $k \geq 2$.
 In view of (\ref{jb}) and (\ref{j'b'}) and the assumptions
$j-b \leq -1$ and $j'-b'\leq -1$, we find that
 $(i'-a')(k+1) \in [1,2k-1]$, where $k \geq 2$.
 So we reach the conclusion that $i'=a'+1$.

 If $j'-b' > -1$, by (\ref{jb}) and (\ref{j'b'}) and the assumptions
$j-b\leq -1$, we find that $-1 < j'-b' \leq k$ and $-k \leq j-b \leq -1$.  It follows that that $(i'-a')(k+1) \in [-k,k-1]$, where $k \geq 1$. Hence we have $i'=a'$. This completes the proof. \qed

It should be noted that the above lemma can be restated
in the form of its converse. Assume that $\pi$ is a permutation in
$\Gamma^k(n,i';j')$. Then $\varphi(\pi)$ is in $\Gamma^k(n,i;j)$.
The verification of this fact is straightforward, and hence it is omitted.

\begin{lem}\label{involution}
 The map $\varphi$ is an involution on $A_{n,k}$, that is, $\varphi^2=I$.
 \end{lem}

\pf
We proceed by induction on $n$. When $n=1$, the lemma is obvious. Suppose that the lemma holds for $n-1$, where $n \geq 2$. We aim to show that it is valid for $n$.

Assume that $\pi$ is a permutation in $A_{n,k}$, that is,
   there exist $0 \leq i \leq n-1$ and $0 \leq j \leq k$ such that $\pi \in \Gamma^k(n,i;j)$. Let $\sigma=\sigma_1 \sigma_2 \cdots \sigma_{n-1}$
    be a permutation of $[n-1]$ such that $\pi=\sigma \leftarrow \pi_n$ with $\pi_n=n-k+j$. Assume $j'$ be the number defined by (\ref{i'}) and (\ref{j'}). Then we have
    \begin{eqnarray*}
 \varphi^2(\pi) &=&\varphi^2(\sigma \leftarrow (n-k+j))\\
 &=& \varphi \big{(}\varphi(\sigma)\leftarrow (n-k+j')\big{)} \\
  &=& \varphi^2(\sigma)\leftarrow (n-k+j) \\
 &=& \sigma\leftarrow (n-k+j)\\
  &=&\pi.                                                           \end{eqnarray*}
This completes the proof.\qed

Combining Lemma \ref{welldefined} and Lemma \ref{involution},
we are led to a bijective proof of Theorem \ref{thmbij},
which is a combinatorial statement of the symmetry property
of $R_{n,k}(x)$ as given by Chung and Graham \cite{ChungGraham}.

\section{Proof of Conjecture \ref{conj-unimodal}}\label{sec:3}

In this section, we give a proof of the conjecture of Hyatt \cite{Hyatt} on the
unimodality of the polynomials $T_k(x)$ associated with
the number of signed permutations with $d$ type $B$ descents and
the type $B$ maximum drop size not larger than $k$.
Based on the polynomials $T_k(x)$, we  define  the polynomials $H_k(x)$ as
given by
\begin{eqnarray}\label{R_k(u)}
\notag
&H_k(x)=&\sum_{l=0}^{k}B_{k-l}(x^{2k+2})(x^{2k+2}-1)^l\sum_{s=l}^{k}{s\choose l}x^{2k+1-s} \\[2pt]
 &&~~+\,\sum_{l=0}^{k}B_{k-l}(x^{-2k-2})(x^{-2k-2}-1)^l\sum_{s=l}^{k}{s\choose l}x^{2(k+1)^2+s}.
\end{eqnarray}

Notice that the sequence of coefficients of $T_k(x)$ is
a subsequence of the sequence of coefficients of $H_k(x)$.
Therefore, the unimodality of $T_k(x)$ follows
from the unimodality of $H_k(x)$.

Let $\widetilde{T}_k(x)=x^k T_k(x)$, that is,
\begin{equation}\label{T'kx}
    \widetilde{T}_k(x)=\sum_{l=0}^{k}B_{k-l}(x^{k+1})(x^{k+1}-1)^l\sum_{i=l}^{k}{i\choose l}x^{k-i}.
\end{equation}

\begin{table}[htp]
\renewcommand{\arraystretch}{1.25}
\begin{tabularx}{130mm}{lX}
\hline
$k$ & $\widetilde{T}_{k}(x)$\\
\hline
0~~ &$1$ \\
1~~&$x+2x^2+x^3$ \\[1pt]
2~~&$x^2+4 x^3+ 6 x^4 + 6 x^5 +4 x^6+2 x^7 +x^8$\\[1pt]
3~~&$x^3+8 x^4+ 12 x^5 + 18 x^6 +23 x^7 + 32 x^8 +32 x^9 +28 x^{10} $$+23 x^{11}$\\
    &~$+\,8 x^{12}+ 4 x^{13} +2 x^{14} + x^{15}$\\
\hline
\end{tabularx}
\caption{The  polynomials $\widetilde{T}_{k}(x)$ for $0\leq k \leq 3$.}
\label{TTk}
\end{table}

The polynomials $\widetilde{T}_k(x)$ for $0\leq k \leq 3$
 are given in Table \ref{TTk}. Like the array representation
 of $Q_k(x)$ given by Chung and Graham \cite{ChungGraham}, we shall use an array representation  of $\widetilde{T}_k(x)$.
The  array representation $t_k$ of the coefficients of $\widetilde{T}_k(x)$
is defined as follows. For $0 \leq i \leq k+1$ and $0 \leq j \leq k$, the $(i,j)$-entry $t_{k}(i,j)$ is
set to be the coefficient of $x^{(k+1)i+j}$ of $\widetilde{T}_k(x)$, that is,
\begin{equation} \label{tk}
    \widetilde{T}_k(x)=\sum_{i=0}^{ k+1} \sum_{j=0 } ^{k} t_{k}(i,j)x^{(k+1)i+j}.
\end{equation}
Similarly, we can arrange the coefficients of $H_k(x)$  in a $ (k+2) \times 2(k+1)$ array $h_k$ so that
\begin{equation*}
    H_k(x)=\sum_{i=0}^{ k+1} ~\sum_{j=0 } ^{2k+1} h_{k}(i,j)x^{2(k+1)i+j}.
\end{equation*}

 It can be seen that the array $h_2$ can be obtained from the array
  $t_2$ by the following operations. First, rotate the array $t_2$ 180 degrees
  counter clockwise. Then put this rotated array in front of $t_2$.
For example,  Table \ref{t2} gives an array $t_2$ and
 Table \ref{r2} gives the corresponding array $h_2$.

\makeatletter\def\@captype{table}\makeatother
\begin{minipage}{.45\textwidth}
\centering
\begin{tabular}{|lcr|}\hline
0 &0  &1\\
4 &6  &6\\
4 &2  &1\\
0 &0  &0\\
\hline
\end{tabular}
\caption{The array $t_2$}
\label{t2}
\end{minipage}
\makeatletter\def\@captype{table}\makeatother
\begin{minipage}{.45\textwidth}
\centering
\begin{tabular}{|lccccr|}\hline
0 &0  &0 &0 &0  &1\\
1 &2  &4 &4 &6  &6\\
6 &6  &4 &4 &2  &1\\
1 &0  &0 &0 &0  &0\\
\hline
\end{tabular}
\caption{The array $h_2$}
\label{r2}
\end{minipage}

In fact, for any $k \geq 0$, $h_k$ can be constructed from $t_k$ in this fashion, which is stated in the following lemma.

\begin{lem}\label{construction}
 For $k \geq 0$, $h_k$  can be obtained by rotating $t_k$  $180$ degrees counter clockwise, and putting the rotated array in front of $t_k$.
\end{lem}

To prove Lemma \ref{construction}, we need the following property.

\begin{lem}\label{insert}
For $k \geq 0$, define
\begin{equation}\label{fk}
F_k(x)=\sum_{l=0}^{k}B_{k-l}(x^{k+2})(x^{k+2}-1)^l\sum_{i=l}^{k}{i\choose l}x^{k+1-i}.
\end{equation}
Arrange the coefficients of $F_k(x)$  in a $ (k+2) \times (k+2)$ array $f_k$ so that
\begin{equation*}
    F_k(x)=\sum_{i=0}^{ k+1} ~\sum_{j=0 } ^{k+1} f_{k}(i,j)x^{(k+2)i+j}.
\end{equation*}
Then the array $f_k$ can be obtained from $t_k$ given in (\ref{tk}) by adding a column of zeros in front of $t_k$.
\end{lem}

\pf
To prove that $f_k$ can be obtained from $t_k$ by inserting a column of
zeros in front of $t_k$, we proceed to verify that
$f_k(i,0)=0$ for $0 \leq i \leq k+1$  and $f_k(i,j+1)=t_k(i,j)$ for $0 \leq i \leq k+1$ and $0 \leq j \leq k$.

For  convenience, for $0 \leq l \leq k$, let
\begin{eqnarray*}
    U_l(t)&=&B_{k-l}(t)(t-1)^l,\\[3pt]
    V_l(t)&=&\sum_{i=l}^{k}{i\choose l}t^{k-i}.
\end{eqnarray*}
Notice that $U_l(t)$ is a polynomial of $t$ of degree $k$ and $V_l(t)$ is
a polynomial of $t$ of degree not larger than $k$.

From the expression (\ref{fk}) of $F_k(x)$, we see that
\[F_k(x)=\sum_{l=0}^{k} xU_l(x^{k+2})V_l(x).\]
Since $U_l(x^{k+2})$  can be seen as a polynomial of $x^{k+2}$ and
the degree of $V_l(x)$ is not larger than $k$, we deduce that
the coefficient of $x^{(k+2)i}$ in $F_k(x)$  equals zero for $0 \leq i \leq k+1$. Hence $f_k(i,0)=0$  for $0 \leq i \leq k+1$.

Next we prove that
 $t_k(i,j)=f_k(i,j+1)$ for $0 \leq i \leq k+1$ and $0 \leq j \leq k$.
We shall adopt the common notation $[x^l]\,p(x)$
for the coefficient of $x^l$ in a polynomial $p(x)$.
It suffices to show that
\begin{equation} \label{coefficient}
[x^{(k+1)i+j}] \,\widetilde{T}_k(x)=[x^{(k+2)i+j+1}]\, F_k(x).
\end{equation}
From the expression (\ref{T'kx}) of $\,\widetilde{T}_k(x)$, it follows that
\[\widetilde{T}_k(x)=\sum_{l=0}^{k} U_l(x^{k+1})V_l(x).\]
Recalling that $V_l(x)$ is a polynomial of $x$ of degree not larger than $k$, for $0 \leq i \leq k+1$ and $0 \leq j \leq k$, it is easily checked that
\begin{eqnarray}\label{x1}
\notag [x^{(k+1)i+j}]\, \widetilde{T}_k(x)&=&\sum_{l=0}^k \Big{(}[x^{(k+1)i}]\,U_l(x^{k+1})\Big{)}
 \Big{(}[x^{j}]\,V_l(x)\Big{)}\\[2pt]
&=&\sum_{l=0}^k \Big{(}[t^i]\,U_l(t)\Big{)}\Big{(}[x^{j}]\,V_l(x)\Big{)}.
\end{eqnarray}

Similarly, we have
\begin{eqnarray}\label{x2}
\notag [x^{(k+2)i+j+1}]\, F_k(x)&=&\sum_{l=0}^k \Big{(}[x^{(k+2)i}]\,U_l(x^{k+2})\Big{)} \Big{(}[x^{j+1}]\,x V_l(x)\Big{)}\\[2pt]
\notag &=&\sum_{l=0}^k \Big{(}[x^{(k+2)i}]\,U_l(x^{k+2})\Big{)} \Big{(}[x^{j}]\,V_l(x)\Big{)}\\[2pt]
&=&\sum_{l=0}^k \Big{(}[t^i]\,U_l(t)\Big{)} \Big{(}[x^{j}]\,V_l(x)\Big{)}.
\end{eqnarray}
Hence  (\ref{coefficient})   follows from (\ref{x1}) and (\ref{x2}).
 So arrive at the conclusion that $f_k(i,j+1)=t_k(i,j)$ for $0 \leq i \leq k+1$ and $0 \leq j \leq k$. This completes the proof. \qed

 We are now ready to give a proof of Lemma \ref{construction}.

\noindent {\it Proof of Lemma \ref{construction}.}
 Write $H_k(x)$ as
 \[ H_k(x)=H'_k(x)+H''_k(x), \] where
\begin{eqnarray}
  H'_k(x) &=& \sum_{l=0}^{k}B_{k-l}(x^{2k+2})(x^{2k+2}-1)^l\sum_{s=l}^{k}{s\choose l}x^{2k+1-s}, \label{H'k} \\[3pt]
  H''_k(x) &=& \sum_{l=0}^{k}B_{k-l}(x^{-2k-2})(x^{-2k-2}-1)^l\sum_{s=l}^{k}{s\choose l}x^{2(k+1)^2+s}.\label{H''k}
\end{eqnarray}
Assume $H'_k(x)$  has the array representation  $h'_k$ such that
\begin{equation*}
    H'_k(x)=\sum_{i=0}^{ k+1} ~\sum_{j=0 } ^{2k+1} h'_{k}(i,j)x^{2(k+1)i+j},
\end{equation*}
and  $H''_k(x)$  has the array representation  $h''_k$ such that
\begin{equation*}
    H''_k(x)=\sum_{i=0}^{ k+1} ~\sum_{j=0 } ^{2k+1} h''_{k}(i,j)x^{2(k+1)i+j}.
\end{equation*}

Clearly, we have $h_k=h'_k+h''_k$.
Using  Lemma \ref{insert} repeatedly, we deduce that
$h'_k$ can be obtained form $t_k$ by adding
$k+1$ columns of zeros in front of $t_k$. Table \ref{hh2} gives an example
of $h'_k$ for $k=2$.

From  the expression (\ref{H'k}) of $H'_k(x)$ and the expression (\ref{H''k}) of $H''_k(x)$, we see that \[ H''_k(x)=H'_k(x^{-1})x^{2(k+1)(k+2)-1}.\]
 Hence, in the form of array representation, we deduce that
 $h''_k$ can be obtained from $h'_k$ by rotating $h'_k$ $180$ degrees.  For example,  the array $h''_2$ in Table \ref{hhh2}
 is constructed from the array $h'_2$ in Table \ref{hh2}.

\makeatletter\def\@captype{table}\makeatother
\begin{minipage}{.45\textwidth}
\centering
\begin{tabular}{|lccccr|}\hline
0 &0 &0 &0 &0  &1\\
0 &0 &0 &4 &6  &6\\
0 &0 &0 &4 &2  &1\\
0 &0 &0 &0 &0  &0\\
\hline
\end{tabular}
\caption{The array $h'_2$}
\label{hh2}
\end{minipage}
\makeatletter\def\@captype{table}\makeatother
\begin{minipage}{.45\textwidth}
\centering
\begin{tabular}{|lccccr|}\hline
0 &0  &0 &0 &0 &0\\
1 &2  &4 &0 &0 &0\\
6 &6  &4 &0 &0 &0\\
1 &0  &0 &0 &0 &0\\
\hline
\end{tabular}
\caption{The array $h''_2$}
\label{hhh2}
\end{minipage}

In view of the fact that $h_k=h'_k +h''_k$ and
the constructions of $h'_k$
and $h''_k$ from $t_k$, we see that the first $k+1$ columns of $h_k$
can be obtained from $t_k$ by a rotation of 180 degrees
and $t_k$ remains to be the last $k+1$ columns of $h_k$.
This completes the proof.\qed

As a consequence of Lemma \ref{construction},
we have  the following property.

\begin{cor}\label{cor}
For $k \geq 0$, the polynomial $H_k(x)$ is symmetric.
\end{cor}

In the array representation,
the symmetry of $H_k(x)$ means that for $0 \leq i \leq k+1$ and $0\leq j \leq 2k+1$, we have
\begin{equation}\label{symmtry}
    h_k(i,j)=h_k(k+1-i,2k+1-j).
\end{equation}

According to Lemma \ref{construction}, it is clear that the sequence of
coefficients of $T_k(x)$ is a subsequence of the coefficients of $H_k(x)$.
We shall prove that  for $k \geq 0$, $H_k(x)$ is a unimodal polynomial.
This implies the unimodality of $T_k(x)$.

\begin{thm}\label{last}
The polynomial $H_k(x)$ is unimodal for any $k \geq 0$.
\end{thm}

 To prove Theorem  \ref{last}, we introduce the polynomials $G_k(x)$
 which will be used to derive a recurrence relation
 of $H_k(x)$.

 Based on the definition (\ref{R_k(u)}) of $H_k(x)$, we define
\begin{eqnarray}\label{G_k(u)}
\notag &G_k(x)=&\frac{1}{x}\sum_{l=0}^{k}B_{k-l}(x^{2k+4})(x^{2k+4}-1)^l\sum_{s=l}^{k}{s\choose l}x^{2k+3-s} \\[2pt]
 &&~~+\sum_{l=0}^{k}B_{k-l}(x^{-2k-4})(x^{-2k-4}-1)^l\sum_{s=l}^{k}{s\choose l}x^{2(k+1)(k+2)+s}.
\end{eqnarray}
Let $g_k$ be the array representation of $G_k(x)$ such that
\begin{equation*}
    G_k(x)=\sum_{i=0}^{ k+1} ~\sum_{j=0 } ^{2k+3} g_{k}(i,j)x^{2(k+2)i+j}.
\end{equation*}
We claim that the array $g_k$ can be obtained from $h_k$ by adding
a column of zeros after the $(k+1)$-st column and adding
 a column of zeros after the $2(k+1)$-st column of $h_k$. The detailed verification of this fact is omitted. Table \ref{h2}
gives an example.

\begin{table}[htp]
\centering
\begin{tabular}{|lccccccr|}\hline
0 &0  &0 &0&0 &0  &1&0\\
1 &2  &4 &0&4 &6  &6&0\\
6 &6  &4 &0&4 &2  &1&0\\
1 &0  &0 &0&0 &0  &0&0\\
\hline
\end{tabular}
\caption{The array $g_2$}
\label{h2}
\end{table}

\begin{lem}\label{recursive}
For   $k \geq 0$, we have
\begin{equation}\label{rec}
    H_{k+1}(x)= G_{k}(x)\cdot ( x +x^2+ \cdots +x^{2k+4})
\end{equation}

\end{lem}

\pf
We aim to  show that
\begin{equation} \label{equg}
    (1-x)\cdot H_{k+1}(x)=x G_{k}(x)\cdot(1-x^{2k+4}),
\end{equation}
which is equivalent to (\ref{rec}).
By the definition of $H_k(x)$ in (\ref{R_k(u)}), we see that  $(1-x)\cdot H_{k+1}(x)$
equals
\begin{eqnarray}\label{ll1}
 \notag\lefteqn{(1-x) \sum_{l=0}^{k+1}B_{k+1-l}(x^{2k+4})(x^{2k+4}-1)^l\sum_{s=l}^{k+1}{s\choose l}x^{2k+3-s} }\\[3pt]
&&+(1-x)\sum_{l=0}^{k+1}B_{k+1-l}(x^{-2k-4})(x^{-2k-4}-1)^l\sum_{s=l}^{k+1}{s\choose l}x^{2(k+2)^2+s}\notag\\[3pt]
   &&~~~=\notag (1-x)\sum_{l=1}^{k+1}B_{k+1-l}(x^{2k+4})(x^{2k+4}-1)^l\sum_{s=l}^{k+1}{s\choose l}x^{2k+3-s} \\[3pt]
   &&~~~~~~~~+(1-x)\sum_{l=1}^{k+1}B_{k+1-l}(x^{-2k-4})(x^{-2k-4}-1)^l\sum_{s=l}^{k+1}{s\choose l}x^{2(k+2)^2+s}\notag\\[3pt]
   &&~~~~~~~~+(1-x)B_{k+1}(x^{2k+4}) \sum_{s=0}^{k+1} x^{2k+3-s} +(1-x)B_{k+1}(x^{-2k-4}) \sum_{s=0}^{k+1} x^{2(k+2)^2+s}\notag\\[3pt]
   &&~~~=-\sum_{l=0}^{k}B_{k-l}(x^{2k+4})(x^{2k+4}-1)^{l+1}\sum_{s=l}^{k}{s\choose {l}}x^{2k+3-s} \notag\\[3pt]
   &&~~~~~~~~+\sum_{l=0}^{k}B_{k-l}(x^{2k+4})(x^{2k+4}-1)^{l+1} {{k+1}\choose {l+1}}x^{k+2} \notag\\[3pt]
    &&~~~~~~~~+\sum_{l=0}^{k}B_{k-l}(x^{-2k-4})(x^{-2k-4}-1)^{l+1} \sum_{s=l}^{k}{s\choose {l}}x^{2(k+2)^2+s+1}\notag\\[3pt]
    &&~~~~~~~~-\sum_{l=0}^{k}B_{k-l}(x^{-2k-4})(x^{-2k-4}-1)^{l+1} {{k+1}\choose {l+1}}x^{(k+2)(2k+5)}\notag\\[5pt]
    &&~~~~~~~~+ B_{k+1}(x^{2k+4})x^{k+2}(1-x^{k+2})
    +B_{k+1}(x^{-2k-4})x^{2(k+2)^2}(1-x^{k+2}).
\end{eqnarray}

On the other hand, by the definition of $G_k(x)$   in  (\ref{G_k(u)}), we
find
\begin{eqnarray*}
x G_{k}(x)\cdot(1-x^{2k+4})&=&
 -\sum_{l=0}^{k}B_{k-l}(x^{2k+4})(x^{2k+4}-1)^{l+1} \sum_{s=l}^{k}{s\choose l}x^{2k+3-s}\\[3pt]
 &&~~~+\sum_{l=0}^{k}B_{k-l}(x^{-2k-4})(x^{-2k-4}-1)^{l+1}\sum_{s=l}^{k}{s\choose l}x^{2(k+2)^2+s+1}.
 \end{eqnarray*}

Comparing the above expression for $x G_{k}(x)\cdot(1-x^{2k+4})$ and
the the first two summations in (\ref{ll1}), to prove (\ref{equg}), it suffices
to show that
the remaining  sum  in (\ref{ll1}) equals zero, that is,
\begin{eqnarray}\label{llll}
\notag\lefteqn{B_{k+1}(x^{2k+4})x^{2k+4}-B_{k+1}(x^{-2k-4})x^{2(k+2)^2}}\\[8pt]
&&=\sum_{l=0}^{k+1}B_{k+1-l}(x^{2k+4})(x^{2k+4}-1)^{l} {{k+1}\choose l}x^{k+2}\notag\\
   &&~~~~-\sum_{l=0}^{k+1}B_{k+1-l}(x^{-2k-4})(x^{-2k-4}-1)^{l} {{k+1}\choose l}x^{(k+2)(2k+5)}.
\end{eqnarray}
It is  known that the type $B$ Eulerian polynomial $B_n(t)$ is a symmetric polynomial of degree
$n$, that is,
\begin{equation*}
    B_{n}(t)=B_{n}(t^{-1})t^n,
\end{equation*} see Brenti \cite{brenti}.
Hence we have
\begin{equation*}
    B_{k+1}(x^{2k+4})x^{2k+4}-B_{k+1}(x^{-2k-4})x^{2(k+2)^2}=0.
\end{equation*}
Thus (\ref{llll}) is equivalent to the following relation
\begin{eqnarray}\label{lll2}
\notag\lefteqn{\sum_{l=0}^{k+1} B_{k+1-l}(x^{2k+4})(x^{2k+4}-1)^{l} {{k+1}\choose l}}\\
     & &=\sum_{l=0}^{k+1}B_{k+1-l}(x^{-2k-4})(x^{-2k-4}-1)^{l} {{k+1}\choose l}x^{2(k+2)^2}.
\end{eqnarray}
Setting $t=x^{2k+4}$ and $n=k+1$, (\ref{lll2}) can be rewritten as
\begin{equation}\label{equ}
    \sum_{l=0}^{n} B_{n-l}(t)(t-1)^{l} {{n}\choose l}=
   \sum_{l=0}^{n} B_{n-l}(t^{-1})(t^{-1}-1)^{l} {{n}\choose l}t^{n+1}.
\end{equation}
To prove (\ref{equ}), we need the following formula
\begin{equation} \label{generating}
    \sum_{n \geq 0} B_n(t) \frac{x^n}{n!}= \frac{(1-t)e^{x(1-t)}}{1-te^{2x(1-t)}},
\end{equation}
which was obtained by Chow and Gessel \cite{Chow}.
Using (\ref{generating}), we get
\begin{eqnarray}\label{las1}
   \notag \lefteqn{ \sum_{n >1} \sum_{j=0}^{n} B_{n-j}(t)(t-1)^j {n \choose j} \frac{x^n}{n!}}\\[3pt]
    &&= \Big{(}\sum_{n \geq 0} B_n(t) \frac{x^n}{n!}\Big{)} \Big{(}\sum_{n \geq 0} (t-1)^n \frac{x^n}{n!}\Big{)}-1 \notag \\[4pt]
    &&=\frac{te^{2x(1-t)}-t}{1-te^{2x(1-t)}}.
\end{eqnarray}
On the other hand, using (\ref{generating}) again, we get
\begin{eqnarray}\label{las2}
   \notag \lefteqn{ \sum_{n >1} \sum_{j=0}^{n} B_{n-j}(t^{-1})(t^{-1}-1)^j {n \choose j} t^{n+1} \frac{x^n}{n!}}\\[3pt]
    &&= t\Big{(}\sum_{n \geq 0} B_n(t^{-1}) \frac{x^n}{n!}\Big{)} \Big{(}\sum_{n \geq 0} (t-1)^n \frac{(tx)^n}{n!}\Big{)}-t \notag\\[4pt]
    &&=\frac{te^{2x(1-t)}-t}{1-te^{2x(1-t)}}.
\end{eqnarray}
Combining (\ref{las1}) and (\ref{las2}), we arrive at  (\ref{equ}). This completes the proof.\qed

Based on Lemma \ref{recursive} and the relationship between
the array representation of $H_k(x)$
and the array representation of $G_k(x)$, we can obtain
the following recurrence relations for the array representation
of $H_k(x)$, which can be verified by induction on $k$. The
detailed proof is omitted.

\begin{cor}\label{corr}
For $0 \leq i \leq k+1$ and $0 \leq j \leq k$, we have
\begin{eqnarray}\label{xxx1}
\notag h_{k}(i,j)&=&h_{k-1}(i,0)+h_{k-1}(i,1)+ \cdots +h_{k-1}(i,j-1)\\[4pt]
&&+h_{k-1}(i-1,j)+h_{k-1}(i-1,j+1)+\cdots+h_{k-1}(i-1,2k-1),
\end{eqnarray}
and for $0 \leq i \leq k+1$ and $ k+1 \leq j \leq 2k+1 $, we have
\begin{eqnarray}\label{xxx2}
\notag h_{k}(i,j)&=&h_{k-1}(i,0)+h_{k-1}(i,1)+ \cdots +h_{k-1}(i,j-2)\\[4pt]
&&+h_{k-1}(i-1,j-1)+h_{k-1}(i-1,j)+\cdots+h_{k-1}(i-1,2k-1),
\end{eqnarray}
where we assume that $h_k(i,j)=0$ when $i <0$.
\end{cor}

Now we are ready to give a proof of Theorem \ref{last}.

\noindent {\it Proof of Theorem \ref{last}.}
We proceed by induction on $k$. By   expression (\ref{R_k(u)}) of $H_k(x)$, we get $H_0(x)=x+x^2$, which is unimodal. Assume that $H_{k-1}(x)$ is unimodal, where $k \geq 1$. We aim to prove that $H_{k}(x)$ is unimodal.

Assume that $k\geq 1$.
Let $(a_0,a_1,\cdots,a_{2k^2+2k-1})$ denote the sequence of coefficients of $H_{k-1}(x)$. By the symmetry of
$H_{k-1}(x)$ as given in Corollary \ref{cor}, we have $a_i=a_{2k^2+2k-1-i}$.
It follows that the unimodality of $H_{k-1}(x)$ is equivalent to the fact that
\begin{equation}\label{z1}
    a_0 \leq a_1 \leq \cdots \leq a_{k^2+k-1}.
\end{equation}
Assume  that $(b_0,b_1,\cdots,b_{2k^2+6k+3})$  is the sequence of
 coefficients of $H_{k}(x)$. By the symmetry of $H_k(x)$, to prove that $H_k(x)$ is unimodal, it suffices for us to prove that
\begin{equation}\label{z2}
    b_0 \leq b_1 \leq \cdots \leq b_{k^2+3k+1}.
\end{equation}

To conduct the induction, we employ  the array representation of $H_k(x)$.
Recall that $h_k$ is the array representation of $H_{k}(x)$ such that
\begin{equation*}
    H_k(x)=\sum_{i=0}^{ k+1} ~\sum_{j=0 } ^{2k+1} h_{k}(i,j)x^{2(k+1)i+j}.
\end{equation*}

Clearly, we have $h_k(i,j)=b_{2(k+1)i+j}$ for $0 \leq i \leq k+1$ and $0 \leq j \leq 2k+1$. Hence we may
 restate (\ref{z2}) in the array representation. More precisely, when $k$ is odd,  (\ref{z2}) can be transformed into the following assertions:
\begin{itemize}
\item[(i)] $h_k(i,j+1)-h_k(i,j) \geq 0$ for $0 \leq i \leq \lfloor \frac{k+2}{2}\rfloor-1 $ and $ 0 \leq j \leq 2k$.

\item[(ii)] $h_k(i,j+1)-h_k(i,j) \geq 0$ for $i=\lfloor \frac{k+2}{2}\rfloor$ and $0\leq j \leq k-1$.

\item[(iii)]
$h_{k}(i,0) - h_{k}(i-1, 2k+1) \geq 0$ for $1 \leq i \leq \lfloor \frac{k+2}{2}\rfloor$.
\end{itemize}
Similarly, when $k$ is even,
 (\ref{z2}) can be recast into the following assertions:
\begin{itemize}
\item[(iv)] $h_k(i,j+1)-h_k(i,j) \geq 0$ for $0 \leq i \leq \frac{k}{2}$ and $ 0 \leq j \leq 2k$.

\item[(v)]
$h_{k}(i,0) - h_{k}(i-1, 2k+1) \geq 0$ for $1 \leq i \leq  \frac{k}{2}$.
\end{itemize}

Now we proceed to prove the above assertions.
It follows from (\ref{xxx1}) that  for $0 \leq i \leq k+1$ and $0 \leq j \leq k-1$,
\begin{equation}\label{z3}
    h_k(i,j+1)-h_k(i,j)=h_{k-1}(i,j)-h_{k-1}(i-1,j).
\end{equation}
Using  (\ref{xxx2}),  we find that for $0 \leq i \leq k+1$ and $k+1 \leq j \leq 2k$,
\begin{equation}\label{z4}
    h_k(i,j+1)-h_k(i,j)=h_{k-1}(i,j-1)-h_{k-1}(i-1,j-1).
\end{equation}
Moreover,  by (\ref{xxx1}) and (\ref{xxx2}), it is easy to check that for $0 \leq i \leq k+1$,
\begin{eqnarray}
h_{k}(i,k) &=& h_{k}(i, k+1)\label{k1}, \\[3pt]
\label{k2} h_{k}(i,0)&=&  h_{k}(i-1, 2k+1).
\end{eqnarray}

We first consider that case when $k$ is odd.  To prove (i),
 we assume that  $\,0 \leq i \leq \lfloor\frac{k+2}{2}\rfloor-1$ and $0 \leq j \leq 2k$.
Here are three subcases. When $0\leq j \leq k-1$,
we claim that $h_k(i,j+1)-h_k(i,j) \geq 0$.
From (\ref{z3}) we see that
\[ h_k(i,j+1)-h_k(i,j)=a_{2ki+j}-a_{2ki-2k+j}.\]
 Since $\,0 \leq i \leq \lfloor\frac{k+2}{2}\rfloor-1$ and $0 \leq j \leq k-1$,
   noting $2 \lfloor \frac{k+2}{2}\rfloor=k+1$, we find that
\[2ki+j \leq 2k\Big{(}\Big{\lfloor}\frac{k+2}{2}\Big{\rfloor}-1\Big{)}+k-1=k^2-1.\]
Clearly, we have $2ki+j \geq 2ki-2k+j$. Thus we may use the
 induction hypothesis to deduce that  $a_{2ki+j}-a_{2ki-2k+j} \geq 0$,
 which is equivalent to the claim.

When  $k+1 \leq j \leq 2k$,
we claim that $h_k(i,j+1)-h_k(i,j) \geq 0$.
By (\ref{z4}),  we get
\[ h_k(i,j+1)-h_k(i,j)=a_{2ki+j-1}-a_{2ki-2k+j-1}.\]
Using the same argument as in the case when $0\leq j\leq k-1$,
 we deduce that
\[2ki+j-1 \leq 2k\Big{(}\Big{\lfloor}\frac{k+2}{2}\Big{\rfloor}-1\Big{)}+2k-1=k^2+k-1.\]
Similarly, we have $2ki+j-1 \geq 2ki-2k+j-1$. Hence  we
 may use the induction hypothesis to deduce that
  $a_{2ki+j-1}-a_{2ki-2k+j-1} \geq 0$, which is equivalent to
   the claim. 
   
   Recall that $h_k(i,k+1)=h_k(i,k)$ for  $\,0 \leq i \leq k+1$ as given in (\ref{k1}).  On the other hand, when $j=k$, assertion (i) becomes the relation $h_k(i,k+1)-h_k(i,k) \geq 0$ for $0 \leq i \leq \lfloor \frac{k+2}{2}\rfloor-1$, which is valid since the equality holds.
Combining the above three cases, assertion (i) is proved.

To prove (ii), we assume that
$i = \lfloor\frac{k+2}{2}\rfloor $ and $0 \leq j \leq k-1$.
We claim that $h_k(i,j+1)-h_k(i,j) \geq 0$.
By (\ref{z3}) and the symmetry relation (\ref{symmtry}), we find that
\begin{eqnarray*}
  h_k(i,j+1)-h_k(i,j) &=& h_{k-1}(i,j)-h_{k-1}(i-1,j) \\[3pt]
   &=& h_{k-1}(k-i,2k-1-j)-h_{k-1}(i-1,j)\\[3pt]
   &=& a_{2k(k-i)+2k-1-j}-a_{2k(i-1)+j}.
\end{eqnarray*}
Since $i = \lfloor\frac{k+2}{2}\rfloor $ and $0 \leq j \leq k-1$, we see that
\[2k(k-i)+2k-1-j \leq 2k\Big{(}k-\Big{\lfloor}\frac{k+2}{2}\Big{\rfloor} \Big{)}+2k-1=k^2+k-1,\]
and
\[2k(k-i)+2k-1-j \geq 2k(i-1)+j.\]
Hence we may use the induction hypothesis to deduce that
$a_{2k(k-i)+2k-1-j}-a_{2k(i-1)+j} \geq 0$. This proves the
claim, and hence assertion (ii) holds.

Note that by (\ref{k2}), we have  $h_{k}(i,0) = h_{k}(i-1, 2k+1)$ for $1 \leq i \leq \lfloor \frac{k+2}{2}\rfloor$. This proves assertion (iii).

Next we turn to the case when $k$ is even.

To prove  (iv), we assume that $0 \leq i \leq \frac{k}{2}$ and $ 0 \leq j \leq 2k$.  When $0 \leq i \leq \frac{k}{2}$ and $0 \leq j \leq k-1$,
we claim that $h_k(i,j+1)-h_k(i,j) \geq 0$.
By (\ref{z3}), we see that
\[ h_k(i,j+1)-h_k(i,j)=a_{2ki+j}-a_{2ki-2k+j}.\]
By the assumptions $\,0 \leq i \leq \frac{k}{2}$ and $0 \leq j \leq k-1$, we see that
\[2ki+j \leq k^2+k-1.\]
Hence we may use the induction hypothesis
to deduce that  $a_{2ki+j}-a_{2ki-2k+j} \geq 0$, which is
equivalent to the claim.

When $0 \leq i \leq \frac{k}{2}-1$ and $k+1 \leq j \leq 2k$,
we claim that $h_k(i,j+1)-h_k(i,j) \geq 0$.
By (\ref{z4}), we find that
\[ h_k(i,j+1)-h_k(i,j)=a_{2ki+j-1}-a_{2ki-2k+j-1}.\]
By the assumptions $\,0 \leq i \leq \frac{k}{2}-1$ and $k+1 \leq j \leq 2k$, we see that
\[2ki+j-1 \leq k^2-1.\]
Hence  the induction hypothesis can be used to get
 $a_{2ki+j-1}-a_{2ki-2k+j-1} \geq 0$, which is equivalent to the claim.

When $i = \frac{k}{2} $ and $k+1 \leq j \leq 2k$,
we claim that $h_k(i,j+1)-h_k(i,j) \geq 0$.
By (\ref{z4}) and the symmetry relation (\ref{symmtry}), we find that
\begin{eqnarray*}
  h_k(i,j+1)-h_k(i,j) &=& h_{k-1}(i,j-1)-h_{k-1}(i-1,j-1) \\[3pt]
   &=& h_{k-1}(k-i,2k-j)-h_{k-1}(i-1,j-1)\\[3pt]
   &=& a_{2k(k-i)+2k-j}-a_{2k(i-1)+j-1}
\end{eqnarray*}
Using the assumptions   $i = \frac{k}{2} $ and $k+1 \leq j \leq 2k$, we get
\[2k(k-i)+2k-j \leq k^2+k-1,\]
and
\[2k(k-i)+2k-j \geq 2k(i-1)+j-1.\]
 Hence  the induction hypothesis can be used to
  deduce that $a_{2k(k-i)+2k-j}-a_{2k(i-1)+j-1} \geq 0$, which
  is equivalent to the claim.

When $j=k$, assertion (iv) takes the form  $h_k(i,k+1)-h_k(i,k) \geq 0$ for $0 \leq i \leq \frac{k}{2}$, which is true
since the equality holds according to (\ref{k1}).
   Combining the above cases, assertion (iv) is proved.
Note that by (\ref{k2}), we have  $h_{k}(i,0) = h_{k}(i-1, 2k+1)$ for $1 \leq i \leq \frac{k}{2}$. This proves assertion (v). So the
proof of the theorem is complete. \qed

\vspace{0.5cm}
 \noindent{\bf Acknowledgments.}  This work was supported by  the 973
Project, the PCSIRT Project of the Ministry of Education,  and the National Science
Foundation of China.

\end{document}